\documentclass[12pt]{article}
\usepackage{geometry, amsmath, amsfonts,diagrams} % see geometry.pdf on how to lay out the page. There's lots.
\diagramstyle[tight,centredisplay%
,dpi=600%              for more modern laser-printers
,noPostScript% 
,heads=LaTeX%
]
\newarrow{Implies} ===={=>}
\newarrow{line} -----
\newtheorem{defi}{Definition}
\newtheorem{lem}{Lemma}
\newtheorem{thm}{Theorem}
\newtheorem{rem}{Remark}

\title{2-gerbes and 2-Tate spaces}
\author{Sergey Arkhipov and Kobi Kremnizer}
\date{} % delete this line to display the current date

\begin{document}

\maketitle
\begin{abstract}
 We construct a central extension of the group of automorphisms of a 2-Tate vector space viewed as a discrete 2-group. This is done using an action of this 2-group on a 2-gerbe of gerbel theories. This central extension is used to define central extensions of double loop groups.
\end{abstract}

\tableofcontents

\section{Introduction}
In this paper we study the question of constructing central extensions of groups using group actions on categories. 

Let $G$ be a group. The basic observation is that the category of $\mathbb{G}_m$ central extensions of
$G$ is equivalent to the category of $\mathbb{G}_m$-gerbes over the classifing stack of $G$. This is in turn equivalent to the category of $\mathbb{G}_m$-gerbes over a point with an action of $G$.
Thus by producing categories with a $G$ action we get central extensions.

We then take this observation one category theoretic level higher. We want to study central extensions of 2-groups. Here a 2-group is a monoidal groupoid such that its set of connected componants is a group with the induced product. We look at the case of a dicrete 2-group, that is we can think of any group $G$ as a 2-group with objects the elements of the group, morphisems the identities and monoidal structure the product. 

We see that $\mathbb{G}_m$-central extensions of a discrete 2-group are the same as 2-gerbes over the classifing stack of the group. This also can be interpreted as a 2-gerbe with $G$ action. Thus to get extensions as a 2-group we should find 2-categories with $G$-action.

These observations are used to define central extensions of automorphism groups of 1-Tate spaces and discrete automorphism 2-groups of 2-Tate spaces. 

The category of $n$-Tate spaces is defined inductively. $0$-Tate spaces are finite dimensional vector spaces. $(n+1)$-Tate spaces are certain indpro objects of the category of $n$-Tate spaces.
To a $1$-Tate space we can associate a 1-gerbe of determinant theories. This 1-gerbe has a natural action of the automorphism group of the 1-tate space. This gives the central extension of the group.

Similarly, to a 2-tate space we can associate a 2-gerbe of gerbel theories with an action of the automorphism group of the 2-Tate space.  This action gives the central extension of the discrete 2-group.

If $G$ is a finite dimensional reductive group and $V$ is a finite dimensional representation we get an embedding of the formal double loop group $G((s))((t))$ into the automorphism group of the 2-tate space $V((s))((t))$. Thus we can restrict the central extension to the double loop group. These central extensions of the double loop group as a 2-group will be used in the future to study the (2-)represntation theory of these groups and relating it to the 2-dimensional Langlands program.

The idea of constructing the higher
central extension in categorical terms belongs
essentially to Michael Kapranov. S.A would like to thank him for sharing the
idea in 2004.

After writing this paper we found out that a similar result was obtained by Osipov in his unpublished
Preprint. S.A. would like to thank Osipov for sharing the manuscript with
him. 

\thanks{K.K. was supported in part by NSF grant DMS-0602007.}
\thanks{S.A. was supported in part by NSERC.} 

\section{Group actions on gerbes and central extensions. }
\subsection{$\mathbb{G}_m$-gerbes and central extensions}
Let`s recall the notion of a group acting on a category.
\begin{defi}\label{groupact}
An action of a group $G$ on a category $C$ consists of a functor $F_g:C\to C$ for each $g\in G$ and a natural transformation $\tau_{g,h}:F_{gh}\to F_g F_h$ s.t.
\begin{equation}\label{action cocycle}
\begin{diagram}
F_{g_1 g_2 g_3} & &\rTo^{\tau_{g_1,g_2 g_3}} & F_{g_1}F_{g_2  g_3} \\ 
  \dTo_{ \tau_{g_1 g_2, g_3}} & & &    \dTo_ {F_{g_1}(\tau_{g_2 ,g_3})} \\
  F_{g_1 g_2}F_{g_3} & &\rTo^{\tau_{g_1 ,g_2}F_{g_3}} & F_{g_1}F_{g_2}F_{g_3}\\ 
\end{diagram}
\end{equation}
commutes for any $g_1,g_2,g_3\in G$.

 We also require that $F_1=Id$ and that $\tau_{1,g}=Id$ and $\tau_{g,1}=Id$.

\end{defi}
Suppose that $C$ is a $\mathbb{G}_m$ gerbe (over a point). By this we mean that:
\begin{itemize}
  \item $C$ is a groupoid. 
  \item $C$ is connected (there exists an arrow between any two objects)
  \item For any object $A$ of $C$, $Aut(A)\simeq \mathbb{G}_m$.
\end{itemize}
Note that this implies that all the $Hom$ spaces are $\mathbb{G}_m$-torsors.

\textit{Remark} If $C$ and $D$ are 1-gerbes then their product $C\times D$ is also a 1-gerbe. This will be used below.

In this case we have the following theorem \cite{Bry}:
\begin{thm}\label{centext}
Let $G$ act on a $\mathbb{G}_m$-gerbe $C$. For each object $A$ of $C$ we get a $\mathbb{G}_m$-central extension $\widetilde{G}_A$. These central extensions depend functorially on $A$ (hence are all isomorphic).
If there exists an equivariant object this extension splits.
\end{thm} 
\textit{Proof:} Let $A\in ob C$. Define 
\begin{equation}
 \widetilde{G}_A=\big\{(g,\phi):g\in G,\phi\in Hom(F_g(A),A)\}
 \end{equation}
  with product  given by
\begin{equation}
(g_1,\phi_1)(g_2,\phi_2)=(g_1 g_2, \phi_1\circ F_{g_1}(\phi_2))
\end{equation}  
Associativity follows from \ref{action cocycle}.

Another way of interpreting this theorem is as follows: An action of $G$ on a gerbe $C$ over a point is the same (by descent) as a gerbe over $\mathbb{BG}$. By taking the cover 
\begin{equation}
\begin{diagram}
pt\\
\dTo\\
\mathbb{BG}\\
\end{diagram}
\end{equation}
we get that such a gerbe gives (again by descent) a line bundle $L$ over $G$ with an isomorphism
\begin{equation}
p_1^*(L)\otimes p_2^*(L) \to m^*(L)
\end{equation}
that gives a product covering that of $G$ and a coherence relation which tells us that this product is associative.
Hence we get
\begin{thm} 
The category of $\mathbb{G}_m$-central extensions of $G$ is equivalent to the category of $\mathbb{G}_m$-gerbes over $\mathbb{BG}$.
\end{thm}

\subsection{Central extension of the automorphism group of a 1-tate space}
Let $\mathcal{V}$ be a 1-Tate space. Recall (or see section \ref{tate}) that this is an ind-pro object in the category of finite dimensional vector spaces, this equivalent to $\mathcal{V}$ having a locally linearly compact topology. Any such is isomorphic to $V((t))$ (formal loops into $V$) but non-canonnically. 
Recall also the notion of a lattice $\mathcal{L}\subseteq \mathcal{V}$ (pro-subspace or linearly compact subspace) and that if $\mathcal{L}_1\subseteq \mathcal{L}_2$ are two lattices then $\mathcal{L}_2 /\mathcal{L}_1$ is finite dimensional.
\begin{defi}\label{det theory}
A determinant theory is  a rule that assigns to each lattice $\mathcal{L}$ a one-dimensional vector space $\Delta_\mathcal{L}$ and to each pair $\mathcal{L}_1\subset\mathcal{L}_2$ an isomorphism
\begin{equation}
\Delta_{\mathcal{L}_1 \mathcal{L}_2}:\Delta_{\mathcal{L}_1}\otimes Det(\mathcal{L}_2/\mathcal{L}_1)\to \Delta_{\mathcal{L}_2}
\end{equation}
such that foe each triple $\mathcal{L}_1\subset\mathcal{L}_2\subset\mathcal{L}_3$ the following diagram commutes
\begin{equation}
\begin{diagram}
\Delta_{\mathcal{L}_1} \otimes Det(\mathcal{L}_2/ \mathcal{L}_1) \otimes Det(\mathcal{L}_3/\mathcal{L}_2) & & & &\rTo & \Delta_{\mathcal{L}_1} \otimes Det(\mathcal{L}_3/ \mathcal{L}_1)\\
\dTo                                                                                                              &                                                                        &    &   &  &  \dTo                                                                                                          \\
\Delta_{\mathcal{L}_2}\otimes Det(\mathcal{L}_3/\mathcal{L}_2) &                                                                          & & & \rTo & \Delta_{\mathcal{L}_3}                                                                          \\
\end{diagram}
\end{equation}
\end{defi}
We have the obvious notion of a morphism between two determinant theories and it is easy to see that the category of determinant theories is in fact a $\mathbb{G}_m$-gerbe.

Let $GL(\mathcal{V})$ be the group of continuous  automorphisms of $\mathcal{V}$. This group acts on the gerbe of determinant theories and hence we get using theorem \ref{centext} a central extension $\widetilde{GL(\mathcal{V})_\mathcal{L}}$
for each choice of lattice $\mathcal{L}$. Unless $\mathcal{V}$ itself is a lattice, this central extension does not split.

\section{Group actions on 2-gerbes and central extensions of 2-groups}
\subsection{2-Groups}
\begin{defi}\label{2group}
A 2-group is a monoidal groupoid $C$ s.t. its set of connected components $\pi_0(C)$ with the induced multiplication is a group.
\end{defi}

The basic example is is the discrete 2-group associated to any group $G$: the set of objects is $G$ itself and morphisms are only the identities. The monoidal structure comes from the group multiplication.
We will denote this discrete 2-group by $\mathcal{G}$.
 
Note that 2-groups can be defined in any category with products (or better in any topos) so we have topological, differential and algebraic 2-groups.

One can define a general notion of extensions of 2-groups but we are only interested in the following case:
\begin{defi}
Let $G$ be a group (in a topos) and $A$ an abelian group (again in the topos). A central extension $\widetilde{\mathcal{G}}$ of the discrete 2-group associated to $G$ by $A$ is a 2-group s.t.:
\begin{itemize}
  \item $\pi_0(\widetilde{\mathcal{G}})\simeq G$
  \item $\pi_1(\widetilde{\mathcal{G}},I)\simeq A$
\end{itemize} 
\end{defi}
Here $I$ is the identity object for the monoidal structure and $\pi_!$ means the automorphism group of the identity object.

\subsection{Action of a group on a bicategory}

Lets recall first the notion of a bicategory (one of the versions of a lax 2-category) \cite{Be}.
\begin{defi}
A bicategory $\mathcal{C}$ is given by:
\begin{itemize}
  \item Objects $A,B,...$
  \item Categories $\mathcal{C}(A,B)$ (whose objects are called 1-arrows and morphisms are called 2-arrows)
  \item Composition functors $\mathcal{C}(A,B)\times\mathcal{C}(B,C)\longrightarrow\mathcal{C}(A,C)$
  \item Natural transformations (associativity constraints)
   \begin{equation}
   \begin{diagram}
   \mathcal{C}(A,B)\times\mathcal{C}(B,C)\times\mathcal{C}(C,D) & & & \rTo & \mathcal{C}(A,B)\times\mathcal{C}(B,D)\\
   &&&&\\
   \dTo & & \rdImplies & &\dTo\\
   &&&&\\
    \mathcal{C}(A,C)\times\mathcal{C}(C,D) & & & \rTo & \mathcal{C}(A,D)\\
   \end{diagram}
   \end{equation}
\end{itemize}
This data should satisfy coherence axioms of the Maclane hexgagon form.
\end{defi}
\textit{Remark} As a bicategory with one object is the same as a monoidal category the coherence axioms should become clear (though long to write).

\begin{defi}
Let $\mathcal{C}$ and $\mathcal{D}$ be two bicategories. A functor $\mathcal{F}:\mathcal{C}\to\mathcal{D}$ consists of:
\begin{itemize}
\item For each object $A\in Ob(\mathcal{C})$ an object $\mathcal{F}(A)\in Ob(\mathcal{D})$
\item A functor $\mathcal{F}_{AB}:\mathcal{C}(A,B)\to\mathcal{D}(\mathcal{F}(A),\mathcal{F}(B))$ for any two objects $A,B\in Ob\mathcal{C}$
\item A natural transformation 
\begin{equation}
\begin{diagram}
\mathcal{C}(A,B)\times\mathcal{C}(B,C) & & & \rTo & \mathcal{C}(A,C)\\
   &&&&\\
   \dTo^{\mathcal{F}_{AB}\times\mathcal{F}_{BC}} & & \rdImplies & &\dTo_{\mathcal{F}_{AC}}\\
   &&&&\\
    \mathcal{D}(\mathcal{F}(A),\mathcal{F}(B))\times\mathcal{D}(\mathcal{F}(B),\mathcal{F}(C)) & & & \rTo & \mathcal{D}(\mathcal{F}(A),\mathcal{F}(C))\\
\end{diagram}
\end{equation}
\end{itemize}
This natural transformation should be compatible with the associativity constraints.
\end{defi}
Again the comparison with monoidal categories should make it clear what are the compatibilities.

\begin{defi}
Let $\mathcal{F}$ and $\mathcal{G}$ be two functors between $\mathcal{C}$ and $\mathcal{D}$. A natural transformation $(\Xi,\xi)$ is given by:
\begin{itemize}
\item A functor $\Xi_{AB}:\mathcal{D}(\mathcal{F}(A),\mathcal{F}(B))\to \mathcal{D}(\mathcal{G}(A),\mathcal{G}(B))$ for each pair of objects
\item A natural transformation 
\begin{equation}
\begin{diagram}
\mathcal{C}(A,B) && & \rTo^{\mathcal{F}_{AB}} & \mathcal{D}(\mathcal{F}(A),\mathcal{F}(B))\\
&\rdline^{\mathcal{G}_{AB}} &&\ldImplies^{\xi_{AB}}& \dTo_{\Xi_{AB}}\\
 &&\rdTo&&\\
 & &&& \mathcal{D}(\mathcal{G}(A),\mathcal{G}(B))  \\
\end{diagram}
\end{equation}
\end{itemize}
These should be compatible with the structures.
\end{defi}
\begin{defi}
Given two natural transformations $(\Xi^1,\xi^1),(\Xi^2,\xi^2):\mathcal{F}\to\mathcal{G}$ a modification is a natural transformation $\phi_{AB}:\Xi^1_{AB}\to\Xi^2_{AB}$ such that
\begin{equation}
\begin{diagram}
\Xi^1_{AB}\mathcal{F}_{AB}&\rImplies^{\xi_{AB}}&\mathcal{G}_{AB}\\
\dImplies^{\phi_{AB}\mathcal{F}_{AB}}&                                             &\dImplies_{Id}\\
\Xi^2_{AB}\mathcal{F}_{AB}&\rImplies^{\xi_{AB}}&\mathcal{G}_{AB}\\
\end{diagram}
\end{equation}
commutes for all $A$ and $B$ and is compatible with all the structures. 
\end{defi}
Now we can define an action of a group on a bicategory:
\begin{defi}
Let $G$ be a group and $\mathcal{C}$ a bicategory. An action of $G$ on $\mathcal{C}$ is given by a functor $\mathcal{F}_g:\mathcal{C}\to\mathcal{C}$ for each $g\in G$ and a 
natural transformation $(\Xi,\xi)_{g,h}:\mathcal{F}_{gh}\to \mathcal{F}_g\mathcal{F}_h$ such that there exists a modification 
\begin{equation}\
\begin{diagram}
\mathcal{F}_{g_1 g_2 g_3} & &\rTo^{(\Xi,\xi)_{g_1,g_2 g_3}} &\mathcal{F}_{g_1}\mathcal{F}_{g_2  g_3} \\  
 \dTo^{ (\Xi,\xi)_{g_1 g_2, g_3}} & & \ldImplies^{\phi_{g_1,g_2,g_3}}&    \dTo_ {\mathcal{F}_{g_1}((\Xi,\xi)_{g_2 ,g_3})} \\
 &&&\\
 \mathcal{F}_{g_1 g_2}\mathcal{F}_{g_3} & &\rTo^{(\Xi,\xi)_{g_1 ,g_2}\mathcal{F}_{g_3}} &\mathcal{F}_{g_1}\mathcal{F}_{g_2}\mathcal{F}_{g_3}\\ 
\end{diagram}
\end{equation}
 for any $g_1,g_2,g_3\in G$ satisfing a cocycle condition.
\end{defi}

\subsection{2-gerbes and central extensions of 2-groups}
Let $A$ be an abelian group.
\begin{defi}
A 2-gerbe (over a point) with band $A$ is a bicategory $\mathcal{C}$ such that
\begin{itemize}
\item It is a 2-groupoid: every 1-arrow is invertible up to a 2-arrow and all 2-arrows are invertible.
\item It is connected: there exists a 1-arrow between any two objects and a 2-arrow between any 1-arrows.
\item The automorphism group of any 1-arrow is isomorphic to $A$
\end{itemize}
\end{defi}
In other words all the categories $\mathcal{C}(A,B)$ are 1-gerbes with band $A$ and the product maps are maps of 1-gerbes. 

\begin{thm}
Suppose $G$ acts on a 2-gerbe $\mathcal{C}$ with band $A$. To this we can associate a central extension $\widetilde{\mathcal{G}}$ of the discrete 2-group associated to $G$ by $A$.
\end{thm}
The construction is the same as in \ref{centext} (with more diagrams to check). A better way of presenting the construction is using descent:
a 2-gerbe with an action of $G$ is the same as a 2-gerbe over $\mathbb{BG}$ (we haven`t defined 2-gerbes in general but the definition is clear \cite{Bre}). Using the same cover as before
$pt\to \mathbb{BG}$ we get a gerbe over $G$ which is multiplicative. That means that we are given an isomorphism
\begin{equation}
p_1^*(\mathcal{F})\otimes p_2^*(\mathcal{F})\to m^*(\mathcal{F})
\end{equation}
satisfying a cocycle condition on the threefold product (here $m:G\times G\to G$ is the multiplication). This gerbe gives in turn an $A$-torsor over $G\times G$ giving the Hom-spaces of the 2-group and the multiplicative structure gives the 
monoidal structure. 
 
This construction also works in the other direction.
Suppose we have a central extension of the discrete 2-group $\mathcal{G}$ associated to the group $G$ by the abelian geoup $A$. 
Then the Hom spaces define an $A$-torsor $\mathcal{HOM}$ over $G\times G$ and the existance of composition means that over $G\times G\times G$ we are given an isomorphism:
\begin{equation}
p_{12}^*(\mathcal{HOM})\otimes p_{23}^*(\mathcal{HOM})\to p_{13}^*(\mathcal{HOM})
\end{equation}
satisfying a cocycle condition over the fourfold product (associativity). Here $p_{ij}$ are the projections. Thus we have a gerbe over $G$ with band $A$. Let`s denote this gerbe by $\mathcal{F}$.

The existence of the monoidal structure implies that we are given an isomorphism over $G\times G$
\begin{equation}
p_1^*(\mathcal{F})\otimes p_2^*(\mathcal{F})\to m^*(\mathcal{F})
\end{equation}
satisfying a cocycle condition on the threefold product. Hence the gerbe is multiplicative. In other words we got:
\begin{lem}
A central extension of the discrete 2-group associated to $G$ by $A$ is the the same as a 2-gerbe over $\mathbb{BG}$ with band $A$. 
\end{lem}
Actually also here we have an equivalence of categories.

\begin{rem} Today`s technology (\cite{Lu})enables one to define n-gerbes with nice descent theory. So we can generalize the whole discussion to:
\begin{thm}\label{n-extension}
The category of n-gerbes with band $A$ and with action of $G$ is equivalent to that of central extensions by $A$ of the discrete n-group associated to $G$.
\end{thm} 
This will be done in another paper.
\end{rem}
\section{2-Tate spaces and 2-groups} \label{tate}
In this section we introduce the notion of a locally compact object introduced by Beilinson and Kato \cite{B,K}.
\subsection{Locally compact objects in a category}
\begin{defi}
Let $C$ be a category. The category of locally compact objects of $C$ is the full subcategory of $Ind(Pro(C))$ consisiting of functors that are isomorphic to diagrams of the following sort: Let $I,J$ be linearly directed orders. Let $F:I^{op}\times J\to C$ be a diagram such that for all $i,i`\in I$ and $j,j`\in J$ $i\leq i`$
and $j\leq j`$ the diagram:
\begin{equation}
\begin{diagram}
F(i`,j)& \rTo & F(i`,j`)\\
\dTo&           & \dTo\\
F(i,j)&\rTo      &F(i`,j)\\
\end{diagram}
\end{equation}
is both cartesian and cocartesian and  vertical arrows are surjections and horizontal arrows are injections. A compact object is a locally compact object isomorphic to one which is constant in the Ind direction. 
\end{defi}

The following statement follows easily from set-theory and the Yoneda lemma:
\begin{lem}
 If $F$ is locally compact then the functors $\stackrel{lim}{\leftarrow}\stackrel{lim}{\rightarrow}F$ and $\stackrel{lim}{\rightarrow}\stackrel{lim}{\leftarrow}F$ are naturally isomorphic.
\end{lem}

From now on we will assume that the indexing sets $I,J$ are countable.

Suppose $C$ is an exact category. Say a sequence of locally compact objects is exact if it can be represented by a map of diagrams $F_1\to F_2\to F_3 :I^{op}\times J\to C$ where all the arrows are exact in $C$. A routine check shows :
\begin{lem}
 The category of locally compact objects of $C$ is exact.
\end{lem}
\begin{rem}
Note that if $C$ is Abelian (and nontrivial) the category of locally compact objects is not abelian.
\end{rem}

Using the standard reindexing trick (Appendix of \cite{AM}) we also get
\begin{lem}
 Let $F_1\to F_2$ be an admissible injection (w.r.t. the exact structure) of compact objects then
 $coker(F_1\to F_2)$ is also a compact object.
\end{lem}
\begin{lem}
 Let $F_1$ and $F_2$ be two admissible compact subobjects of $F$ then $F_1\times_F F_2$ is also compact. 
\end{lem}
Now we can define inductively n-Tate spaces (we still assume that the indexing sets are countable):
\begin{defi}
A 0-Tate space is a finite dimensional vector space. Suppose we have defined the category of n-Tate spaces. A $n+1$-Tate space is a locally compact object of n-Tate spaces. A lattice of a $n+1$-Tate space is an admissible compact subobject.
\end{defi}
Note that any 2-Tate space is of the form $\mathcal{V}((t))$ where $\mathcal{V}$ is a 1-Tate space. An example of a lattice in this case is $\mathcal{V}[[t]]$.

\subsection{Some facts on 1-Tate spaces}
 We have from the previous section that:
\begin{lem} 
The category of (1-)Tate spaces is an exact category with injections set-theoretic injections and surjections dense morphisms.
\end{lem}
Recall also the notion of the determinant grebe associated to a Tate space $\mathcal{V}$. From now on we will denote it by $\mathcal{D}_\mathcal{V}$.
\begin{lem}
Let
\begin{equation}
0\to \mathcal{V}'\to \mathcal{V} \to \mathcal{V}'' \to 0
\end{equation} 
be an admissible exact sequence of Tate spaces. Then we have an equivalence of $\mathbb{G}_m$-gerbes
\begin{equation} 
  \mathcal{D}_{\mathcal{V}'}\otimes \mathcal{D}_{\mathcal{V}''}\to \mathcal{D}_\mathcal{V}
\end{equation}
such that if $\mathcal{V}_1\subset\mathcal{V}_2\subset\mathcal{V}_3$ then we have a natural transformation
\begin{equation}
\begin{diagram}
\mathcal{D}_{\mathcal{V}_1} \otimes \mathcal{D}_{\mathcal{V}_2/ \mathcal{V}_1} \otimes \mathcal{D}_{\mathcal{V}_3/\mathcal{V}_2} & &  &\rTo & \mathcal{D}_{\mathcal{V}_1} \otimes \mathcal{D}_{\mathcal{V}_3/ \mathcal{V}_1}\\
\dTo                                                                                                              &                                                                           & &  &  \dTo                                                                                                          \\
&&\ldImplies&&\\
&&&&\\
\mathcal{D}_{\mathcal{V}_2}\otimes \mathcal{D}_{\mathcal{V}_3/\mathcal{V}_2} &                                                                           & & \rTo & \mathcal{D}_{\mathcal{V}_3}                                                                          \\
\end{diagram}
\end{equation}
and if we have $\mathcal{V}_1\subset\mathcal{V}_2\subset{V}_3\subset\mathcal{V}_4$ then a cubical diagram of natural transformations commutes.
\end{lem}     
\subsection{2-Tate spaces and gerbel theories}

It follows from the previous discussion that:
\begin{lem}
Let $\mathbb{V}$ be a 2-Tate space.
\begin{enumerate}
  \item If $\mathbb{L}'\subset \mathbb{L}$ are two lattices then $\mathbb{L}/\mathbb{L}'$ is a 1-Tate space. 
  \item For any two lattices $\mathbb{L}$ and $\mathbb{L}'$ there exists a third lattice $\mathbb{L}''\subset \mathbb{L}\cap\mathbb{L}'$
\end{enumerate}  
\end{lem}
Now we can define a gerbel theory.
\begin{defi}
Let $\mathbb{V}$ be a 2-vector space. A gerbel theory $\mathbb{D}$ is 
\begin{itemize}
\item For each lattice $\mathbb{L}\subset\mathbb{V}$ a $\mathbb{G}_m$-gerbe $\mathbb{D}_\mathbb{L}$
\item If $\mathbb{L}'\subset\mathbb{L}$ are two lattices then we have an equivalence
\begin{equation}
\begin{diagram}
\mathbb{D}_\mathbb{L} &\rTo^{\phi_{\mathbb{L}\mathbb{L}'}}& &\mathbb{D}_{\mathbb{L}'}\otimes \mathcal{D}_{\mathbb{L}/\mathbb{L}'}\\
\end{diagram}
\end{equation}
\item  For $\mathcal{V}_1\subset\mathcal{V}_2\subset\mathcal{V}_3$ we have a natural transformation
\begin{equation}
\begin{diagram}
\mathbb{D}_{\mathcal{V}_1} \otimes \mathcal{D}_{\mathcal{V}_2/ \mathcal{V}_1} \otimes \mathcal{D}_{\mathcal{V}_3/\mathcal{V}_2} & &  &\rTo & \mathbb{D}_{\mathcal{V}_1} \otimes \mathcal{D}_{\mathcal{V}_3/ \mathcal{V}_1}\\
\dTo                                                                                                              &                                                                           & &  &  \dTo                                                                                                          \\
&&\ldImplies&&\\
&&&&\\
\mathbb{D}_{\mathcal{V}_2}\otimes \mathcal{D}_{\mathcal{V}_3/\mathcal{V}_2} &                                                                           & & \rTo & \mathbb{D}_{\mathcal{V}_3}                                                                          \\
\end{diagram}
\end{equation}
\end{itemize}
Given $\mathcal{V}_1\subset\mathcal{V}_2\subset\mathcal{V}_3$ these natural transformations should commute on a cubical diagram.
\end{defi}

Now we have
\begin{thm}
Gerbel theories on a given 2-Tate space $\mathbb{V}$  form a $\mathbb{G}_m$ 2-gerbe $\mathbb{GERB}_\mathbb{V}$. 
\end{thm}
Let`s denote $\mathbb{GL}(\mathbb{V})$ the group of continuos automorphisms of a  2-Tate space $\mathbb{V}$. This group acts naturally on the 2-gerbe  $\mathbb{GERB}_\mathbb{V}$. 
\textit{Remark} the action is actually a strict one. 
We get:
\begin{thm}
Let $\mathbb{V}$ be a 2-Tate space. Given a lattice $\mathbb{L}\subset\mathbb{V}$ we get a $\mathbb{G}_m$ central extension of the discrete 2-group associated to $\mathbb{GL}(\mathbb{V})$.
\end{thm}

\begin{rem}
 Using \ref{n-extension} we can go on and define central extensions of discrete n-groups of automorphism of n-Tate spaces. 
\end{rem}

\subsubsection{Application: central extension of a double loop group}
Let $G$ be a finite dimensional reductive group over a field. Let $V$ be a finite dimensinal representation of $G$. From this data we get a map 
\begin{equation}
G((s))((t))\to \mathbb{GL}(V((s))((t)))
\end{equation}
 where $G((s))((t))$ is the formal double loop group of $G$. From this embedding we get a central extension of the discrete 2-group $G((s))((t))$.
\subsubsection{ A variant}
There is another way to think about $\mathbb{G}_m$-gerbes. 
\begin{defi}
 Let $\mathit{Pic}$ be the symmetric monoidal groupoid of 1-dimensional vector spaces. A $\mathbb{G}_m$-gerbe is a module category over this monoidal category equivalent to $\mathit{Pic}$ as module categories (where $\mathit{Pic}$ acts on itself by the monoidal structure). 
\end{defi}
  
This definition is equivalent to the definition given before. Now, following Drinfeld \cite{Dr} we define a graded version of a $\mathbb{G}_m$-gerbe.

\begin{defi}
 Let $\mathit{Pic}^\mathbb{Z}$ be the symmetric monidal groupoid of $\mathbb{Z}$-graded 1-dimensional vector spaces with the super-commutativity constraint ($a\otimes b \to (-1)^{deg(a)deg(b)} b\otimes a$). A $\mathbb{Z}$-graded $\mathbb{G}_m$-gerbe is a module category over  $\mathit{Pic}^\mathbb{Z}$ equivalent to it as module categories.
\end{defi}

We have a map from $\mathit{Pic}^\mathbb{Z}$ to the discrete 2-group $\mathbb{Z}$ which sends a 1-dimensional graded vector space to its degree. This map induces a functor between graded $\mathbb{G}_m$-gerbes and $\mathbb{Z}$-torsors. We can now repeat the entire story with $\mathbb{Z}$-graded gerebs. For instance, instead of a determinant theory we will get a graded determinant theory. The $\mathbb{Z}$-torsor coing to it will be the well known dimension torsor of dimension theories. A dimension theory for a 1-Tate space is a rule of associating an integer to each lattice satisfying similar conditions as a determinant theory.

In this way we will get for a 2-Tate space an action of $\mathbb{GL}(\mathbb{V})$ on the $\mathbb{G}_m$-gerbe of dimension torsors. This action will give us a central extension of the group $\mathbb{GL}(\mathbb{V})$ (not the 2-group!). And similarly we can get central extensions of groups of the form $G((s))((t))$. Thus we see that if we work with graded determinant theory we get a central extension of the dicrete 2-group $\mathbb{GL}(\mathbb{V})$ which induces the central extension of the group $\mathbb{GL}(\mathbb{V})$ (For this central extension see \cite{O}).

\begin{rem}
 Another reason to work with graded theories is that they behave much better for the direct sum of 1-Tate spaces. It is true that the determinant gerbe of the direct sum of 1-Tate spaces is equivalent to the tensor product of the gerbes but this equivalence depends on the ordering. If one works with graded determinant theories this equivalence will be canonical. 
\end{rem}

\end{document}